\documentclass[11pt]{amsart}
\usepackage{geometry}                
\usepackage{csquotes}
\usepackage{graphicx}
\usepackage{amssymb}
\usepackage{amsthm}
\usepackage{epstopdf}
\theoremstyle{plain}

\newtheorem{example}{Example}
\newtheorem{rem}{Remark}

\title[McKean--Pontryagin for optimal transport]{A McKean--Pontrygin maximum principle for entropic-regularized optimal transport}
\author{Sebastian Reich}
\address{Institut f\"ur Mathematik, Universit\"at Potsdam}
\date{March 2026}

\begin{document}

\maketitle

\begin{abstract}
    This note outlines a mean-field approach to dynamic optimal transport problems based on the recently proposed McKean--Pontryagin maximum principle. Key aspects of the proposed methodology include i) avoidance of sampling over stochastic paths, ii) a fully variational approach leading to constrained Hamiltonian equations of motion, and iii) a unified treatment of deterministic and stochastic optimal transport problems. We also discuss connections to well-known dynamic formulations in terms of forward-backward stochastic differential equations and extensions beyond classical entropic-regularized transport problems.
\end{abstract}

%
\section{Problem statement}
%

We wish to bridge two given distributions $\pi_0$ and $\pi_T$ in a variable $x \in \mathbb{R}^{d_x}$ along the controlled stochastic differential equation
\begin{equation} \label{eq:SDE}
{\rm d}\tilde X_t = U_t{\rm d}t + \sqrt{2}\Sigma^{1/2}{\rm d}\tilde B_t, 
\end{equation}
subject to minimizing the cost
\begin{equation} \label{eq:cost U}
\mathcal{J}(U) = \frac{1}{2} \int_0^T \mathbb{E}_{\tilde X_t}\|U_t \|^2_R \,{\rm d}t
\end{equation}
in the control $U = U_{[0,T]}$, $U_t \in \mathbb{R}^{d_u}$, $d_u=d_x$, and  under constraints $\tilde X_0\sim \pi_0$ and $\tilde X_t\sim \pi_T$. Here, $\|u\|_R$ denotes the weighted norm defined by
\begin{equation}
    \|u\|^2_R = u^{\rm T}R^{-1}u = R^{-1}: u u^{\rm T}
\end{equation}
for given symmetric positive definite matrix $R$. The problem reduces to the classical optimal transport problem for $\Sigma = 0$ and $R=I$ and to the Schr\"odinger bridge problem for $\Sigma = R=I$ \cite{peyré2019computational,CGP21,Tang26}.

The goal of this note is to formulate a mean-field formulation of the Schr\"odinger bridge problem based on an extension of the classical Pontryagin maximum principle \cite{pontryagin}. More precisely, the proposed formulation extends the McKean--Pontryagin maximum principle, first introduced in \cite{R25b}, to the general coupling problem (\ref{eq:SDE})--(\ref{eq:cost U}).

We note that the subsequent formulations can be easily extended to stochastic differential equations of the form
\begin{equation} \label{eq:general SDE}
    {\rm d}\tilde X_t = b(\tilde X_t){\rm d}t + GU_t{\rm d}t +
    \sqrt{2} \Sigma^{1/2}{\rm d}\tilde B_t
\end{equation}
for given $b(x) \in \mathbb{R}^{d_x}$ and $G\in \mathbb{R}^{d_x\times d_u}$. More specifically, consider controlled underdampled Langevin dynamics, which leads to $X_t = (Q_t^{\rm T},V_t^{\rm T})^{\rm T}$, 
\begin{equation}
    b(q,v) = \left( \begin{array}{c} v\\ 0 \end{array}\right), \quad
    G = \left(\begin{array}{c} 0\\ I \end{array}\right), \quad 
    \Sigma = \left( \begin{array}{cc} 0 & 0 \\ 0 & \tilde \sigma I \end{array}\right).
\end{equation}

One can also extend the variational formulation to multiplicative noise $\Sigma(x)$ and mean-field formulations such as the Kalman--Wasserstein gradient flow \cite{CRS22}, which corresponds to $R = \Sigma = C(\rho_t)$ with $C(\rho_t)$ the covariance matrix of $X_t\sim \rho_t$. 

%
\section{McKean--Pontryagin maximum principle}
%

We apply the McKean--Pontryagin variational formulation from \cite{R25b} with states $X_t(a) \in \mathbb{R}^{d_x}$ and co-states $P_t(a) \in \mathbb{R}^{d_x}$. The distribution of $X_t(a)$ is denoted by $\rho_t$. Here, $a \in \mathbb{R}^{d_x}$ are labels that remain constant and serve as independent variables. The proposed variational formulation will ensure that $X_t$ and $\tilde X_t$ agree in law for all $t\in [0,T]$.

In a first step, we introduce initial and final conditions in a weak form; that is,
\begin{equation}
    X_0(a) = b_0, \qquad X_T(a) = b_T
\end{equation}
for suitable boundary values $b_0 = b_0(a)$ and $b_T = b_T(a)$ such that
\begin{equation}
    \int_{\mathbb{R}^{d_x}} f(X_0(a)) \rho_0(a){\rm d}a = \int_{\mathbb{R}^{d_x}} f(x)\pi_0(x){\rm d}x = \pi_0[f]
\end{equation}
and
\begin{equation}
    \int_{\mathbb{R}^{d_x}} f(X_T(a)) \rho_0(a){\rm d}a = \int_{\mathbb{R}^{d_x}} f(x)\pi_T(x){\rm d}x = \pi_T[f]
\end{equation}
for all suitable test functions $f(x)$. This formulation takes into account that there is a relabeling symmetry, which preserves the prescribed marginal densities $\pi_0$ and $\pi_T$.

Next, following \cite{R25b}, we state the action functional $\mathcal{S} = \mathcal{S}(X,P,U,\beta,\psi)$ defined by
\begin{subequations} \label{eq:action}
    \begin{align}
\mathcal{S} =& \,\int_0^T \int_{\mathbb{R}^{d_x}}
\left( \langle P_t,\dot{X}_t-U_t\rangle - \Sigma : D_x^2 \psi_t(X_t))
+ \frac{1}{2} \|U_t\|_R^2  \right) \rho_0(a){\rm d}a \,{\rm d}t\\
&  \,-\,\int_0^T \int_{\mathbb{R}^{d_x}}\langle P_t-\nabla_x \psi_t(X_t),\beta_t \rangle \,\rho_0(a){\rm d}a
\,{\rm d}t \\
&\,-\, \int_{\mathbb{R}^{d_x}} \psi_T(X_T) \,\rho_0(a){\rm d}a 
+\pi_T[\psi_T] +\int_{\mathbb{R}^{d_x}} \psi_0(X_0)\,\rho_0(a){\rm d}a 
- \pi_0[\psi_0]
    \end{align}
\end{subequations}
with inner product $\langle a,b\rangle =
a^{\rm T} b$. We introduce the Hamiltonian $\mathcal{H} = \mathcal{H}(X,P,U,\beta,\psi)$
\begin{equation} \label{eq:Hamiltonian}
    \mathcal{H} = \int_{\mathbb{R}^{d_x}}\left(
    \langle P,U \rangle + \Sigma : D_x^2 \psi(X) - \frac{1}{2} \|U\|^2_R
    + \langle \beta,P-\nabla_x \psi(X) \rangle \right) \,\rho_0(a) {\rm d}a
\end{equation}
and note that
\begin{subequations}
\begin{align}
\mathcal{S} =&\, \int_0^T \left\{
\int_{\mathbb{R}^{d_x}}
\langle P_t,\dot{X}_t\rangle \,\rho_0(a)
{\rm d}a - \mathcal{H}(X_t,P_t,U_t,\beta_t,\psi_t)\right\} {\rm d}t\\
& \,-\, \int_{\mathbb{R}^{d_x}} \psi_T(X_T) \,\rho_0(a){\rm d}a 
+\pi_T[\psi_T] +\int_{\mathbb{R}^{d_x}} \psi_0(X_0)\,\rho_0(a){\rm d}a 
- \pi_0[\psi_0].
\end{align}
\end{subequations}
We note that $\Sigma : D_x^2 \psi(x)$ in (\ref{eq:Hamiltonian}) acts as a state-dependent running cost similar to stochastic optimal control \cite{Carmona}; the main difference being that the cost is defined implicitly through the relation $P_t = \nabla_x \psi_t(X_t)$. As will be verified in the following, $\psi_t(x)$ solves the associated Hamilton--Jacobi--Bellman equation \cite{Carmona}.

Taking variations of (\ref{eq:action}) with respect to $P_t$, $X_t$, $\beta_t$, $\psi_t$, and $U_t$, the evolution equations
\begin{subequations} \label{eq:evolution equations}
    \begin{align}
        \dot{X}_t &= \frac{\delta \mathcal{H}}{\delta P_t} = U_t + \beta_t,\\
        -\dot{P}_t &= \frac{\delta \mathcal{H}}{\delta X_t} = \nabla_x (\Sigma : D_x^2 \psi_t(X_t)) -
        D_x^2 \psi_t(X_t) \,\beta_t
    \end{align}
\end{subequations}
emerge subject to the following constraints:
\begin{subequations} \label{eq:constraints}
    \begin{align}
        0&= \frac{\delta \mathcal{H}}{\delta \beta_t} =P_t - \nabla_x \psi_t(X_t),\\
        0 &= \frac{\delta \mathcal{H}}{\delta \psi_t} = \nabla_x \cdot \left(\rho_t(X_t)\left\{\beta_t +
        \Sigma \nabla_x \log \rho_t(X_t)\right\}\right),\\
        0 &= \frac{\delta \mathcal{H}}{\delta U_t} =P_t - R^{-1}U_t.
    \end{align}
\end{subequations}
Here, it has been used in (\ref{eq:constraints}b) that
\begin{subequations}
\begin{align}
\int_{\mathbb{R}^{d_x}} \Sigma : D_x^2 \psi_t(X_t(a))\,\rho_0(a){\rm d}a &= 
\int_{\mathbb{R}^{d_x}} \Sigma : D_x^2 \psi_t(x)\,\rho_t(x) {\rm d}x \\ &= 
\int \nabla_x \cdot (\rho_t(x)\, \Sigma \nabla_x \log \rho_t(x))\,\psi_t(x)\,{\rm d}x
\end{align}
\end{subequations}
and
\begin{equation}
\int_{\mathbb{R}^{d_x}} \langle \beta_t(X_t(a)),\nabla_x \psi_t(X_t(a))\rangle \,\rho_0(a){\rm d}a =
- \int_{\mathbb{R}^{d_x}} \nabla_x \cdot(\rho_t(x) \,\beta_t(x)) \,\psi_t(x)\,{\rm d}x.
\end{equation}
The boundary conditions are given by
\begin{equation}
    \rho_0 = \pi_0, \qquad \rho_T = \pi_T, \qquad P_0 = \nabla_x \psi_0, \qquad
    P_T = \nabla_x \psi_T,
\end{equation}
which arise from variations with respect to $X_0$, $\psi_0$ and $X_T$, $\psi_T$, respectively. Here we have assumed for simplicity that $\rho_0(x) = \pi_0(x)>0$ everywhere. The potentials $\psi_0$ and $\psi_T$ are unknown and arise from the boundary constraints on the density $\rho_t$. We also find that $\beta_t$ is not uniquely determined by (\ref{eq:constraints}b); but a natural choice is
\begin{equation}
\beta_t = -\Sigma \nabla_x \log \rho_t(X_t).
\end{equation}
The gauge freedom in the choice of $\beta_t$ is a consequence of the relabelling symmetry of our McKean--Pontryagin formulation. We will demonstrate below that, in fact, $\beta_t$ can be arbitrarily chosen while still providing the optimal control law $u_t(x) = R\nabla_x \psi_t(x)$ \cite{R25b}.

It follows from the Lagrangian variational principle \cite{holmI} that the Hamiltonian (\ref{eq:Hamiltonian}) is preserved and takes the constant value
\begin{equation} \label{eq:energy}
    \mathcal{E} =
    \int_{\mathbb{R}^{d}}\left(   \frac{1}{2}\|U_t\|_R^2 + 
    \Sigma : D_x^2 \psi_t(X_t)
     \right) \rho_0(a)\,{\rm d}a
\end{equation}
along solutions of (\ref{eq:evolution equations}).

The optimal transport formulation of Benamou--Brenier is recovered for $\Sigma \to 0$ and $R = I$ \cite{Tang26}. Furthermore, for $\Sigma = R = I$, we have found a (deterministic) mean-field formulation of the dynamic Schr\"odinger bridge problem. Indeed, the law, $\rho_t(x)$, of $X_t$ satisfies the Fokker--Planck equation
\begin{equation} \label{eq:Liouville}
\partial_t \rho_t = \mathcal{L}_{u_t}^\dagger \rho_t,
\end{equation}
while $\psi_t(x)$ satisfies the Hamilton--Jacobi--Bellman equation
\begin{equation} \label{eq:HJB}
-\partial_t \psi_t = \mathcal{L}_{u_t}\psi_t - \frac{1}{2}\|u_t\|_R^2
\end{equation}
with generator $\mathcal{L}_{u_t}f = \langle \nabla_x f,u_t \rangle + \Sigma: D_x^2 f$ and control law $u_t(x) = R\nabla_x \psi_t(x)$ \cite{Carmona}. More precisely, it follows from
\begin{equation}
\dot{P}_t = \frac{{\rm d}}{{\rm d}t}\nabla_x \psi_t(X_t) = 
\partial_t \nabla_x \psi_t(X_t) + D_x^2\psi_t(X_t)\dot{X}_t
\end{equation}
together with the evolution equations (\ref{eq:evolution equations}) for $X_t$ and $P_t$ that
$\nabla_x \psi_t$ satisfies
\begin{equation} \label{eq:gradient HJB}
    -\partial_t \nabla_x \psi_t(X_t) = D_x^2 \psi_t(X_t)U_t + \nabla_x (\Sigma: D_x^2 \psi_t(X_t)),
\end{equation}
which implies (\ref{eq:HJB}). Furthermore, this calculation demonstrates that the time evolution of $\psi_t(x)$ is independent of the specific choice of $\beta_t$ as already observed in \cite{R25b}.

Note that the potentials $\psi_t$ are directly related to the Schr\"odinger potentials $\phi_t$ and $\hat \phi_t$ \cite{Tang26} via
\begin{equation}
    \hat \phi_t = e^{\psi_t}, \qquad \rho_t = \phi_t \hat \phi_t.
\end{equation}
for all $t \in [0,T]$.

\begin{rem}
    It is possible to replace the potential $\psi_t$ and $\nabla_x \psi_t$ in the variational principle by a vector-valued function $\Psi_t$ leading to the modified Hamiltonian
    \begin{equation}  \label{eq:Hamiltonian_2}
    \mathcal{H} = \int_{\mathbb{R}^{d_x}}\left(
    \langle P,U \rangle + \Sigma : D_x \Psi(X) - \frac{1}{2} \|U\|^2_R
    + \langle \beta,P-\Psi(X) \rangle \right)\, \rho_0(a) {\rm d}a.
    \end{equation}
Taking variations with respect to $\Psi_t$ then leads to the constraint
\begin{equation}
    0 = \beta_t + \Sigma \nabla_x \log \rho_t(X_t).
\end{equation}
Furthermore, a regularizing penality term 
\begin{equation} \label{eq:penality}
    \mathcal{J}(\Psi) = \frac{\gamma}{2} \int_{\mathbb{R}^{d_x}}\|P - \Psi(X)\|^2\, \rho_0(a){\rm d}a,
\end{equation}
$\gamma>0$, can be added to the Hamiltonian (\ref{eq:Hamiltonian_2}) without changing the resulting equations of motion \cite{Optmization}. In this context, we briefly describe a connection to adjoint matching \cite{AM,Tang26}. Since $U_t = R\,\Psi_t(X_t)$, the constraint (\ref{eq:constraints}a) or, equivalently, minimizing (\ref{eq:penality})
can be viewed as matching the adjoint $P_t$ and the control $U_t$. We again mention the observation from \cite{R25b} that the resulting control law $u_t(x) = R\, \Psi_t(x)$ is independent of the specific choice of $\beta_t$ as long as $\Psi_T$ is specified appropriately. Finally, if the assumption is dropped that $\Psi_t$ is itself the gradient of a potential $\psi_t$, ideal barotropic fluid dynamics \cite{Salmon88} can also be treated within this framework by setting $\Sigma \equiv 0$ and adding the internal energy $e(\rho_t)$ as a cost function to the Hamiltonian (\ref{eq:Hamiltonian_2}). It would be of interest to discuss the effect of a non-zero $\Sigma$ in the context of barotropic fluid dynamics.
\end{rem}

\begin{example}
Consider the specific case of a given target distribution $\pi_0$, which one wishes to transform into a standard Gaussian at time $T>0$; that is, $\pi_T = {\rm N}(0,I)$. Then, one can consider the following iterative procedure to find the optimal control $U_t$, which follows the idea of iterative proportional fitting \cite{Tang26}. A natural initial guess is $U_t = 0$,
which formally follows from $P_t = 0$ for $t\in [0,T]$. The resulting distribution of $X_t$ at time $T$ is denoted by $\rho_T^{(0)}$, suggesting the following corrected boundary condition
\begin{equation}
\psi_T^{(1)}(x)= - \frac{1}{2}\| x\|^2 -
\log \rho_T^{(0)}(x) = 
\log \pi_T(x) -
\log \rho_T^{(0)}(x)
\end{equation}
and $P_T(a) = \nabla \psi_T^{(1)}(X_T(a))$.
Upon induction, one derives
\begin{equation}
\psi_T^{(K)}(x) = K\,\log \pi_T(x ) - \sum_{k=0}^{K-1}
\log \rho_T^{(k)}(x).
\end{equation}
and $P_T(a)= \nabla \psi_T^{(K)}(X_T(a))$. The process converges as $\nabla_x \log \rho_T^{(k)} \to \nabla_x \log \pi_T$.
\end{example}

\section{A non-variational Lagrangian approach}
The following closely related approach has been proposed in \cite{ZOL24}. In the setting of this note, the formulation is given by
\begin{subequations} \label{eq:Lagrangian HJB}
    \begin{align}
        \dot{X}_t =&\, U_t - \Sigma \nabla_x \log \rho_t(X_t),\\
        -\dot{Y}_t =&\, -\frac{1}{2}\|U_t\|^2_R + \Sigma: D_x^2 \psi_t(X_t)
        + \langle \nabla_x \psi_t(X_t),\Sigma \nabla_x \log \rho_t (X_t)\rangle,\\
        0 =& \,Y_t - \psi_t(X_t),\\
        0 =& \,R^{-1} U_t - \nabla_x \psi_t(X_t). 
    \end{align}
\end{subequations}
Since
\begin{equation}
\dot{Y}_t = \partial_t \psi_t(X_t) + \langle \nabla_x \psi_t(X_t),\dot{X}_t\rangle,
\end{equation}
equation (\ref{eq:Lagrangian HJB}b) constitutes a Lagrangian formulation of (\ref{eq:HJB}) along solutions of (\ref{eq:Lagrangian HJB}a). However, the formulation does not satisfy a Hamiltonian variational principle \cite{holmI}.

\section{Stochastic McKean--Pontryagin maximum principle}
The time evolution of $\psi_t(x)$ does not depend on the specific choice of $\beta_t$. This fact has already been exploited in the numerical implementations proposed in \cite{R25b}. Here we use this freedom to build a link to forward-backward SDE formulations \cite{Carmona} of stochastic optimal control in the context of the Schr\"odinger bridge problem. Hence, let us formally set 
\begin{equation}
\beta_t = \sqrt{2}\Sigma^{1/2}\dot{B}_t
\end{equation}
with $B_t$ denoting $d$-dimensional Brownian motion. Stochastic integrals are interpreted in Stratonovitch form \cite{Pavliotis2016}. Then (\ref{eq:evolution equations})-(\ref{eq:constraints}) give rise to the pathwise formulation
\begin{subequations} \label{eq:FBSDE}
    \begin{align}
        {\rm d}X_t &= \,U_t\,{\rm d}t + {\rm d}\beta_t,\\
        -{\rm d}P_t &=\, \nabla_x (\Sigma : D_x^2 \psi_t(X_t))\,{\rm d}t -
        D_x^2 \psi_t(X_t) \circ {\rm d}\beta_t,\\
        0&= \,P_t - \nabla_x \psi_t(X_t),\\
        0 &= \,P_t - R^{-1}U_t
    \end{align}
\end{subequations}
with $X_t = X_t(a,B_{[0,T]})$ and $P_t = P_t(a,B_{[0,T]})$. However, while delivering the optimal control law $u_t(x) = R \nabla_x \psi_t(x)$, this stochastic formulation is not optimal from a variational perspective.

Furthermore, replacing the Stratonovitch integral in (\ref{eq:FBSDE}) by its It\^o interpretation, one obtains a set of evolution equations formally equivalent to the forward-backward stochastic differential equations associated with stochastic optimal control \cite{Carmona}; that is,
\begin{subequations} 
    \begin{align}
        {\rm d}X_t &= \,U_t\,{\rm d}t + {\rm d}\beta_t,\\
        -{\rm d}P_t &=\,-
        D_x^2 \psi_t(X_t) \, {\rm d}\beta_t,\\
        0&= \,P_t - \nabla_x \psi_t(X_t),\\
        0 &= \,P_t - R^{-1}U_t.
    \end{align}
\end{subequations}
It is important to note that $\psi_t(x)$ is a deterministic function despite the fact that both $P_t$ and $X_t$ depend on $B_{[0,T]}$. A related formal equivalence has also been observed in the context of (\ref{eq:Lagrangian HJB}) 
\cite{ZOL24}.

\section{Closing remarks}

An extension to stochastic differential equations (\ref{eq:general SDE}) subject to given marginal distributions fits naturally within the proposed framework and provides an alternative to approaches based on forward-backward stochastic differential equations \cite{Carmona}. It is possible to consider diffusion matrices $\Sigma$ and drift terms $b$, which depend on the law, $\rho_t$, of $X_t$.  

It would also be of interest to treat diffusion matrix $\Sigma$ as an additional control term with, for example, reward 
\begin{equation}
-\frac{\gamma}{2} \|\Sigma\|^2_{\rm F} = -\frac{\gamma}{2} \Sigma : \Sigma,
\end{equation}
$\gamma>0$. A formal application of Pontryagin's maximum principle to an appropriately extended Hamiltonian (\ref{eq:Hamiltonian}) then leads to the optimal choice
\begin{equation}
\Sigma(x) = \frac{1}{\gamma} D_x^2 \psi_x(x)
\end{equation}
and (\ref{eq:energy}) becomes
\begin{subequations}
\begin{align}
     \mathcal{E} &=
    \int_{\mathbb{R}^{d}}\left(   \frac{1}{2}\|U_t\|_R^2 + \frac{1}{2\gamma}
    D_x^2 \psi_t(X_t) : D_x^2 \psi_t(X_t)
     \right) \rho_0\,{\rm d}a\\
    &= \int_{\mathbb{R}^{d}}\left(   \frac{1}{2}\|R\nabla_x\psi_t(x)\|_R^2 + \frac{1}{2\gamma}
    \|D_x^2 \psi_t(x)\|_{\rm F}^2
     \right) \rho_t(x)\,{\rm d}x,
\end{align}
\end{subequations}
which places a smoothness penalty on $\psi_t$. Note that our mean-field approach does not require the symmetric matrix $\Sigma$ to be positive definite; however, the connection to a stochastic differential equation of the form (\ref{eq:SDE}) is no longer straightforward. 

Discrete-time interacting particle approximations can be obtained by discretizing the action (\ref{eq:action}) in label space and time \cite{Runge-Kutta}. See also \cite{R25b} for numerical approximations of $\beta_t$ and $\psi_t$. The mean-field formulation (\ref{eq:evolution equations})-(\ref{eq:constraints}), based on (\ref{eq:action}), is expected to lead to less noisy approximations than its stochastic counterpart (\ref{eq:FBSDE}). A detailed numerical comparison of different methodologies based on mean-field ordinary differential equations and forward-backward stochastic differential equations for optimal transport problems is left for future work.

\medskip

\paragraph{Acknowledgments.}

This work has been partially funded by Deutsche Forschungsgemeinschaft (DFG) - Project-ID 318763901 - SFB1294. 

\bibliographystyle{unsrt}
\bibliography{bib-database}

\end{document}